\documentclass[11pt]{article}
\usepackage{amsmath, amscd, amsthm, amssymb, graphics, mathrsfs,textcomp}
\usepackage{fancyhdr}
\usepackage[pagebackref=true]{hyperref}
\hypersetup{backref}

\newcommand{\setof}[1]{\{ #1 \}}

\newcommand{\aone}{{\mathbb A}^1}
\newcommand{\pone}{{\mathbb P}^1}

\renewcommand{\O}{{\mathcal O}}

\newcommand{\F}{{\mathcal F}}

\newcommand{\cplx}{{\mathbb C}}

\newcommand{\Z}{{\mathbb Z}}

\newcommand{\ga}{{\mathbb G}_{\bf a}}
\newcommand{\gm}{{\mathbb G}_{\bf m}}

\newcommand{\isomto}{\stackrel{\sim}{\longrightarrow}}

\newcommand{\Spec}{\operatorname{Spec}}

\newcommand{\Sm}{{\mathcal Sm}}

\theoremstyle{plain}
\newtheorem{thm}{Theorem}[section]
\newtheorem{lem}[thm]{Lemma}

\newtheorem{cor}[thm]{Corollary}
\newtheorem{prop}[thm]{Proposition}

\theoremstyle{definition}

\theoremstyle{remark}
\newtheorem{rem}[thm]{Remark}

\numberwithin{equation}{section}

\newcommand{\shrinkmargins}[1]{
  \addtolength{\textheight}{#1\topmargin}
  \addtolength{\textheight}{#1\topmargin}
  \addtolength{\textwidth}{#1\oddsidemargin}
  \addtolength{\textwidth}{#1\evensidemargin}
  \addtolength{\topmargin}{-#1\topmargin}
  \addtolength{\oddsidemargin}{-#1\oddsidemargin}
  \addtolength{\evensidemargin}{-#1\evensidemargin}
  }

\shrinkmargins{.9}

\begin{document}
\pagestyle{fancy}
\renewcommand{\sectionmark}[1]{\markright{\thesection\ #1}}
\fancyhead{}
\fancyhead[LO,RE]{\bfseries\footnotesize\thepage}
\fancyhead[LE]{\bfseries\footnotesize\rightmark}
\fancyhead[RO]{\bfseries\footnotesize\rightmark}
\chead[]{}
\cfoot[]{}
\setlength{\headheight}{1cm}

\title{{\bf Vector bundles on contractible smooth schemes}}
\author{Aravind Asok \\ \begin{footnotesize}Department of Mathematics\end{footnotesize} \\ \begin{footnotesize}University of Washington\end{footnotesize} \\ \begin{footnotesize}Seattle, WA 98195 \end{footnotesize} \\ \begin{footnotesize}\url{asok@math.washington.edu}\end{footnotesize}\\
\and Brent Doran\footnote{This material is based upon work supported
by the National Science Foundation, agreement No. DMS-0111298.}
\\ \begin{footnotesize}School of Mathematics\end{footnotesize} \\
\begin{footnotesize}Institute
for Advanced Study\end{footnotesize} \\
\begin{footnotesize}Princeton,
NJ 08540\end{footnotesize} \\
\begin{footnotesize}\url{doranb@math.ias.edu}\end{footnotesize}}

\date{}

\maketitle

\begin{abstract}
We discuss algebraic vector bundles on smooth $k$-schemes $X$
contractible from the standpoint of $\aone$-homotopy theory; when $k
= \cplx$, the smooth manifolds $X(\cplx)$ are contractible as
topological spaces. The integral algebraic K-theory and integral
motivic cohomology of such schemes are that of $\Spec k$.  One might
hope that furthermore, and in analogy with the classification of
topological vector bundles on manifolds, algebraic vector bundles on
such schemes are all isomorphic to trivial bundles; this is almost
certainly true when the scheme is affine. However, in the non-affine
case this is false: we show that (essentially) every smooth
$\aone$-contractible strictly quasi-affine scheme that admits a
$U$-torsor whose total space is affine, for $U$ a unipotent group,
possesses a non-trivial vector bundle.  Indeed we produce explicit
arbitrary dimensional families of non-isomorphic such schemes, with
each scheme in the family equipped with ``as many" (i.e., arbitrary
dimensional moduli of) non-isomorphic vector bundles, of every
sufficiently large rank $n$, as one desires; neither the schemes nor
the vector bundles on them are distinguishable by algebraic
K-theory.  We also discuss the triviality of vector bundles for
certain smooth complex affine varieties whose underlying complex
manifolds are contractible, but that are not necessarily
$\aone$-contractible.
\end{abstract}

\section{Introduction}
In this note, we study the set of isomorphism classes of vector
bundles on smooth $k$-schemes that are contractible in the sense of
$\aone$-homotopy theory (as introduced in \cite{MV}); such schemes
will be called $\aone$-contractible.  We wish to stress three
counter-intuitive points. First, as our main results show, there are
lots of these, both of the schemes and of the bundles on a typical
fixed such scheme (see Theorem \ref{thm:main}, and Corollaries
\ref{cor:discrete} and \ref{cor:moduli}). Second, they arise quite
naturally and explicitly, so should not be considered pathological.
Third, the standard cohomology theories (at least those theories
representable on the $\aone$-homotopy category) are completely
insensitive to these structures, and so are missing a surprising
amount of algebro-geometric data.

Regarding the third point let us be more specific right from the
start. Since motivic cohomology is representable in the
$\aone$-homotopy category (see \cite{VDelNotes} Theorem
2.3.1\footnote{The proof of this fact requires, at the moment, that
$k$ be a perfect field.}), $\aone$-contractible schemes have the
motivic cohomology of $\Spec k$ and so, for instance, have no
non-trivial algebraic cycles. Similarly, and more importantly for
our present purposes, since algebraic K-theory is representable in
the $\aone$-homotopy category (see \cite{MV} \S 4 Theorem 3.13), one
knows that the algebraic K-theory of any $\aone$-contractible smooth
$k$-scheme is isomorphic to that of $\Spec k$; already from
$K_0(\Spec k) \cong \Z$ this implies that all vector bundles are
stably trivial.

Given an $\aone$-contractible smooth scheme $X$, it is therefore
natural to ask whether all the vector bundles on $X$ are in fact
trivial, especially given that topological vector bundles on open
contractible manifolds are trivial. Indeed, recalling the
Quillen-Suslin theorem for affine space (itself the prototypical
smooth $\aone$-contractible scheme, and the only one known before
\cite{AD1}), one may view this as a generalized Serre problem. We
show there is a stark dichotomy between the affine and strictly
quasi-affine cases: in the affine case, the answer seems to be yes,
whereas in the quasi-affine case we prove that the answer is a
resounding no and construct explicit counter-examples in abundance.
What is especially interesting is that none of the standard means
for distinguishing vector bundles on a scheme (e.g., Chern classes,
algebraic K-theory, algebraic cycles) can play any role at all; from
their standpoint, all the bundles are indistinguishable from a
trivial bundle.

On the other hand, it follows from Corollary
\ref{cor:strictqafftrivvb} that there are moduli of strictly
quasi-affine surfaces which are not $\aone$-contractible (nor are
the complex surfaces even contractible in the sense of manifolds)
and yet they admit only trivial vector bundles.  Thus having
non-trivial vector bundles is by no means a necessary feature of
being strictly quasi-affine.

\subsubsection*{Representability properties of the functor ``isomorphism classes of vector bundles"}
We put the above discussion in a broader context.  Let $\Sm_k$
denote the category of separated, finite type, smooth schemes
defined over $k$.  The $\aone$-homotopy category is constructed by
embedding the category $\Sm_k$ in a larger category of {\em spaces},
equipping that category with the structure of a model category and
then forming the associated homotopy category.  The category of
spaces is taken to be the category of simplicial Nisnevich sheaves
on $\Sm_k$. The homotopy category can be formed by localizing along
two classes of morphisms: first, along the {\em simplicial weak
equivalences} and second along the {\em $\aone$-weak equivalences}.
We refer the reader to (\cite{MV} \S 2 Theorem 3.2) for precise
details regarding this construction.  Here and through the remainder
of the paper $[\cdot,\cdot]_s$ and $[\cdot,\cdot]_{\aone}$ will
denote the set of simplicial homotopy classes of maps and
$\aone$-homotopy classes of maps between spaces.

Let ${\mathscr V}(\cdot)$ (resp. ${\mathscr V}_n(\cdot)$) denote the
functor that assigns to an object $X \in \Sm_k$ the set of
isomorphism classes of (rank $n$) locally free sheaves on $X$.  Let
$BGL_n$ denote the usual simplicial classifying space defined in
\cite{MV} \S 4.1, then by {\em ibid.} \S 4 Proposition 1.16, one
knows that the set of simplicial homotopy classes of maps
$[X,BGL_n]_{s}$ can be identified with $H^1_{Nis}(X,GL_n)$.  Using a
version of ``Hilbert's Theorem 90" (i.e., that $GL_n$ is a ``special
group" in the sense of Serre), we know that the last group is
isomorphic to $H^1_{Zar}(X,GL_n)$ which is, essentially by
construction, isomorphic to ${\mathscr V}_n(X)$.  Ideally, one hopes
that ${\mathscr V}_n(X)$ descends to a functor on the
$\aone$-homotopy category and is representable by the space $BGL_n$,
i.e., that $[X,BGL_n]_{\aone} = {\mathscr V}_n(X)$.

\subsubsection*{Positive results}

In the case $n = 1$ this ideal scenario is the reality, without
restriction on $X$.  Recall that a smooth scheme $X$ is called
$\aone$-rigid (see \cite{MV} \S 3 Example 2.4), if for any smooth
scheme $U$, the map $Hom_{\Sm_k}(U,X) \longrightarrow Hom_{\Sm_k}(U
\times \aone,X)$ induced by pullback along the projection $U \times
\aone \longrightarrow U$ is a bijection.   Morel and Voevodsky show
that, since the sheaf $\gm$ is $\aone$-rigid, $B\gm = BGL_1$ is in
fact $\aone$-local (see {\em ibid.} \S 3 Definition 2.1).  Thus,
$[X,B\gm]_{\aone} = [X,B\gm]_{s}$ and one concludes ${\mathscr
V}_1(X) = H^1_{Zar}(X, \gm) = [X,B\gm]_{\aone}$ (see {\em ibid.} \S 4
Proposition 3.8).

Furthermore, Morel argues (see \cite{MBundle} Theorem 3) that if one
restricts ${\mathscr V}_n(\cdot)$ to a functor on the category of
smooth {\em affine} schemes then again this ideal is realized, at
least if $n \neq 2$ (it is expected that $n = 2$ works as well, but
the details remain to be written out).  Consequently, if $X$ is an
affine $\aone$-contractible smooth $k$-scheme, then every vector
bundle on $X$ (of rank $n \neq 2$) is isomorphic to a trivial bundle.

\subsubsection*{Negative results}

Unfortunately, for $n \geq 2$, the functors ${\mathscr V}_n(X)$
cannot descend to functors on the homotopy category without a
restriction on $X$: it has long been known (and was pointed out to
us by Morel) that even with $X = \pone$, the canonical map
${\mathscr V}(\pone) \longrightarrow {\mathscr V}(\pone \times
\aone)$ induced by pull-back via the projection morphism is {\em
not} a bijection, as we discuss in \S \ref{ss:qproj}.  Observe that
this means the space $BGL_n$ is {\em not} $\aone$-local for $n > 1$
({\em cf.}, \cite{MV} p. 138).  Indeed, one can show that
$\aone$-locality of $BGL_n$ is equivalent to the assertion that, for
any smooth scheme $X$, and $i = 0,1$, the canonical map
$H^i_{Nis}(X,GL_n) \longrightarrow H^i_{Nis}(X \times \aone,GL_n)$
is a bijection (combine Proposition 1.16 of \cite{MV} \S 4 and
\cite{MIntro} Lemma 3.2.1); Morel has called this latter
cohomological condition on a group ``strong $\aone$-invariance."

Nevertheless, Morel's results might make one hope that some form of
homotopy invariance holds for the functor ${\mathscr V}_n(X)$ for
general $n$ beyond the affine case.  For instance, perhaps any
$\aone$-weak equivalence of smooth schemes $f: X \rightarrow Y$
where $X$ is affine would induce a bijection $f^*: \mathscr{V}(Y)
\rightarrow \mathscr{V}(X)$; when, in addition, $Y$ is affine this
is true by the discussion above.

\begin{rem}
Indeed, by the Jouanolou-Thomason homotopy lemma (see e.g.,
\cite{Weibel} Proposition 4.4), given any smooth scheme $Y$
admitting an ample family of line bundles (e.g., a quasi-projective
variety), there exists a smooth affine scheme $X$ and a Zariski locally
trivial smooth morphism with fibers isomorphic to affine spaces $f:
X \longrightarrow Y$.  In particular, this morphism is an
$\aone$-weak equivalence, so the above na\"ive hope would reduce the
study of vector bundles on such schemes to the case of affine
varieties!  Unfortunately, Theorem \ref{thm:main} shows that this is false.

Nevertheless, Morel's results combined with the Jouanolou-Thomason
homotopy lemma give the following general picture. Suppose $Y$ is a
smooth scheme over a perfect field $k$ (according to our
conventions, this means $Y$ is separated, regular and Noetherian and
thus admits an ample family of line bundles).  As long as $n$ is a
strictly positive integer $\neq 2$, then for any smooth affine
scheme $X$ that is $\aone$-weakly equivalent to $Y$, we have a
bijection $[Y,BGL_n]_{\aone} \cong {\mathscr V}_n(X)$.
\end{rem}


Alternatively, homotopy invariance might hold for a slightly broader
class of varieties than affine ones, say for quasi-affine schemes
with ``nice enough" affine closures.  In \cite{AD1}, using
techniques for studying unipotent group actions developed in
\cite{DK}, we constructed many examples in characteristic $0$ of
non-isomorphic strictly quasi-affine $\aone$-contractible smooth
schemes.  Using this construction for arbitrary $k$, we will see
that neither of the above generalizations are possible; it seems
Morel's results are in fact the strongest one can expect.
Furthermore, we will see that any attempt to quantify the lack of
homotopy invariance must account for arbitrarily many non-isomorphic
vector bundles. Specifically, our main goal in this paper is to
prove the following result.

\begin{thm}
\label{thm:main} Let $k$ be a field.  Suppose $X$ is a finite type,
smooth, affine $\aone$-contractible $k$-scheme equipped with a free
everywhere stable action of a split connected unipotent group $U$.
\begin{itemize}
\item[i)] The quotient $X/U$ exists as a smooth $\aone$-contractible
quasi-affine scheme.
\item[ii)] If $X/U$ is affine, then for every positive integer $n$, the pull-back map ${\mathscr V}_n(X/U) \longrightarrow {\mathscr V}_n(X)$ is a bijection.
\item[iii)] If $X$ is isomorphic to affine space and $X/U$ is affine, then every vector bundle on $X/U$ is isomorphic to a trivial bundle.
\item[iv)] If $X/U$ is not affine, but admits a smooth quasi-affine closure with at least one codimension $\geq 2$ boundary component, then $X/U$ admits non-trivial vector bundles of rank $m$ for all sufficiently large $m$.
\end{itemize}
\end{thm}

\begin{rem}
One might suspect that any quasi-affine scheme that is not affine
has non-trivial vector bundles, but this is false in general.
Indeed, one can show that all vector bundles on the complement of
finitely many points in ${\mathbb A}^2$ are trivial. (We will
generalize this fact in Corollary \ref{cor:strictqafftrivvb}.)
\end{rem}

\begin{rem}
The boundary component condition in $(iv)$ above is imposed for ease
of proof; it almost certainly can be removed, and is satisfied for
instance when $\Spec k[X]^U$ is smooth, which is true in the generic
case of our construction.  We have asked (see \cite{AD1,AD2}), in
analogy with the structure theory of contractible manifolds, whether
any smooth $\aone$-contractible variety can be realized as such a
quotient of affine space by the free action of a unipotent group.  A
positive solution to this question would have the following
consequence: removing the boundary component condition in Theorem
\ref{thm:main} implies the interesting dichotomy that {\em all}
non-affine smooth $\aone$-contractible varieties have a non-trivial
vector bundle, whereas {\em all} affine ones would only have trivial
vector bundles.
\end{rem}

In \S \ref{s:examples}, we will expand on this theorem by placing a
lower bound on ``how many" non-trivial vector bundles such a
quasi-affine $\aone$-contractible variety can have, and thereby
construct large dimensional families of examples with as many
non-trivial vector bundles as one likes, all indistinguishable from
the trivial bundle from the point of view of algebraic K-theory (or
of any invariant representable in the $\aone$-homotopy category).

\subsubsection*{Contractible complex affine algebraic varieties}
In addition,  we will visit the generalized Serre problem as
discussed in \cite{Za} \S 8.  We will say that a scheme $X$ over
$\cplx$ is topologically contractible if $X(\cplx)$ equipped with
its usual structure of a complex manifold is contractible as a
topological space.  The generalized Serre problem asks: if $X$ is a
smooth complex affine algebraic variety which is topologically
contractible, then are all algebraic vector bundles on $X$
isomorphic to trivial bundles?  Note that if $X$ is an
$\aone$-contractible smooth scheme over $\cplx$, then $X$ is
necessarily topologically contractible (see \cite{AD1} Lemma 2.5).
However, not all topologically contractible complex varieties are
$\aone$-contractible (see \cite{AD2}); for example any topologically
contractible smooth complex surface of log-general type is {\em not}
$\aone$-contractible (in fact such surfaces can be shown to be {\em
$\aone$-rigid}).  We will observe in \S 4, putting together results
of several authors, that the generalized Serre problem is true for
all topologically contractible smooth complex varieties of dimension
$\leq 2$; consequently there are positive dimensional moduli of
smooth surfaces that each admit only trivial vector bundles.
Finally, we will present some examples of topologically contractible
smooth complex $3$-folds all of whose vector bundles are isomorphic
to trivial bundles.

\subsubsection*{Conventions and Definitions}
The word ``scheme" will mean separated scheme, locally of finite
type over a field $k$.   The word ``variety" will mean reduced,
finite type scheme.  A scheme $X$ is called $\aone$-contractible if
the canonical morphism $X \longrightarrow \Spec k$ is an
$\aone$-weak equivalence in the sense of \cite{MV} \S 3 Definition
2.1.  A scheme $X$ is called quasi-affine if there exists an affine
scheme $\bar{X}$ and an open immersion $X \hookrightarrow \bar{X}$;
we will refer to quasi-affine schemes that are not affine as {\em
strictly quasi-affine} schemes.  If $X$ is any scheme, we let ${\sf
Vec}(X)$ denote the category of finite rank locally free
$\O_X$-modules and, as above, ${\mathscr V}(X)$ will denote the set
of isomorphism classes of vector bundles on $X$.

Throughout, $U$ will denote a split {\em connected} unipotent
$k$-group.  Splitness of $U$ implies that $U$ admits an increasing
filtration by normal subgroups with sub-quotients isomorphic to
$\ga$, and in particular that $U$ is isomorphic to affine space as a
$k$-scheme.  Observe that if $\operatorname{char}(k) > 0$, then
split unipotent groups {\em can} have non-trivial finite subgroups
(e.g., the kernel of the Artin-Schreier morphism $\ga
\longrightarrow \ga$).

Actions of groups on schemes are always assumed to be left actions;
actions will be called free if they are scheme-theoretically free,
i.e., the action morphism is a closed immersion.  If $X$ is a scheme
equipped with an action of $U$, then $X/U$ will denote the geometric
quotient of $U$ by $X$, if it exists as a scheme.

A {\em $U$-torsor} over a scheme $X$ will be a triple $({\mathscr
P},\pi,U)$ consisting of a faithfully flat, finite presentation
morphism $\pi: {\mathscr P} \longrightarrow X$, from a left
$U$-scheme ${\mathscr P}$, such that the canonical morphism $U
\times {\mathscr P} \longrightarrow {\mathscr P} \times {\mathscr
P}$ is an isomorphism onto ${\mathscr P} \times_X {\mathscr P}$.
Observe that in this situation, $U$ acts freely on ${\mathscr P}$
(see \cite{GIT} Lemma 0.6) and $X$ is a geometric quotient of
${\mathscr P}$ by $U$.

Our notation and terminology will follow \cite{AD1} unless otherwise
mentioned. However, the reader need not be familiar with the results
of {\em ibid.}, as long as she takes on faith Theorems 3.10 and
4.11 therein: in essence, (a) there is a (computable) notion of an
everywhere stable $U$-action on an affine scheme $X$, (b) it is
equivalent to $X$ being endowed with the structure of a $U$-torsor
over a quasi-affine scheme $X/U$, and (c) in certain circumstances
we can explicitly identify the complement of the open immersion of
$X/U$ in $\Spec k[X/U]$ using geometric invariant theory.

\subsubsection*{Acknowledgements}

We would like to thank Fabien Morel for interesting discussions
around the topic of vector bundles in $\aone$-homotopy theory, and
in particular for pointing out the quasi-projective counter-example
for rank $2$ vector bundles. We would also like to thank Jacob Lurie
for a useful conversation that helped simplify our proofs. Both of
these interactions were facilitated by a workshop in Topology at the
Banff International Research Station.  Finally, we would like to
thank the referee for useful comments.

\section{Vector bundles and $U$-torsors}
If $q: X \rightarrow X/U$ is a $U$-torsor, we observe
the induced map $q^*: \mathscr{V}(X/U) \rightarrow \mathscr{V}(X)$
can have very different character depending on whether $X/U$ is
affine or non-affine. The affine and strictly quasi-affine cases
will be used to prove Theorem \ref{thm:main}.

\subsubsection*{The affine case: $q^*$ is a bijection}
\begin{lem}
\label{lem:affinequot}
Suppose $q: X \rightarrow X/U$ is a $U$-torsor with
$X/U$ a smooth, affine scheme.  Then $q$ induces a bijection
$$
q^*: {\mathscr V}(X/U) \isomto {\mathscr V}(X).
$$
If in addition $X$ is isomorphic to affine space, then every vector
bundle on $X/U$ is isomorphic to a trivial bundle.
\end{lem}

\begin{proof}
Lindel proved (see \cite{Lindel}) that if $Y$ is a smooth affine
$k$-scheme then pullback via the projection map $Y \times {\mathbb
A}^n \longrightarrow Y$ induces a bijection ${\mathscr V}(Y) \isomto
{\mathscr V}(Y \times {\mathbb A}^n)$.

According to the hypotheses, $q: X \longrightarrow X/U$ equips the
triple $(X,q,U)$ with the structure of a $U$-torsor over $X/U$ and
$X/U$ is affine.  Observe that for any affine scheme $Y$,
$H^1(Y,\ga) = H^1(Y,\O_Y) = 0$ by \cite{EGAIII.1} Th\'eor\`eme
1.3.1.  As $U$ is split, an inductive argument shows that $H^1(Y,U)
= 0$ for any such $Y$. Thus $(X,q,U)$ must be a trivial $U$-torsor
over $X/U$, whence $X \cong U \times X/U$.  Thus, the first result
follows from the discussion of the previous paragraph.\footnote{See
the proof of Corollary 3.2 in the Appendix to \cite{AD1} for details
(note that because of the splitness assumption, there is no
restriction on the base field).}

The final statement, where $X$ is assumed to be affine space, now
follows from the Quillen-Suslin theorem (see e.g., \cite{Quillen})
that all vector bundles on affine space are isomorphic to trivial
bundles.
\end{proof}

\subsubsection*{A quasi-projective counterexample: failure of surjectivity}

\label{ss:qproj} Consider the projection morphism $p_1: \pone \times
\aone \longrightarrow \pone$.  We will show that the pull-back map
$p_1^*:{\mathscr V}(\pone) \longrightarrow {\mathscr V}(\pone \times
\aone)$ is not a bijection. By Grothendieck's description of the
category of vector bundles on $\pone$, we know that every locally
free sheaf on $\pone$ is isomorphic to a direct sum of rank $1$
locally free sheaves.  A vector bundle on $\pone \times \aone$
isomorphic to a pull-back of a vector bundle $\F$ on $\pone$ is
necessarily isomorphic to the (external) tensor product of $\F$ and
$\O_{\aone}$.  There is a rank $2$ vector bundle ${\mathcal E}$ on
$\pone \times \aone$ whose restriction to $\pone \times \setof{0}$
is trivial and whose restriction to $\pone \times \setof{1}$ is
isomorphic to $\O(1) \oplus \O(-1)$.  This means ${\mathcal E}$ is
not isomorphic to the pull-back of any bundle on $\pone$.

\subsubsection*{The strictly quasi-affine case: failure of injectivity}
Now assume $q: X \longrightarrow X/U$ is a $U$-torsor with $X/U$
{\em not} affine.  The existence of $\aone$-contractible strictly quasi-affine $X/U$ will be proved in \S \ref{s:examples}.  More generally, for the rest of this subsection, and in particular for the statements of Lemmas \ref{lem:short}, \ref{lem:exist} and \ref{lem:dense}, we assume we are in the following situation:
\begin{itemize}
\item[i)] $X/U$ is an open dense subscheme of a finite type smooth scheme $\overline{X/U}$, with the inclusion denoted $j: X/U \hookrightarrow \overline{X/U}$,
\item[ii)] we denote by $Z$ the closed complement of $X/U$ in $\overline{X/U}$ equipped with the reduced induced scheme structure and assume it is {\em non-empty}.
\end{itemize}
In this situation, we have a localization sequence in $G$-theory (see \cite{Sr} Proposition 5.15):
$$
\cdots \longrightarrow G_1(X/U) \longrightarrow G_0(Z)
\longrightarrow G_0(\overline{X/U}) \longrightarrow G_0(X/U)
\longrightarrow 0.
$$
As both $X/U$ and $\overline{X/U}$ are finite type smooth schemes,
we know by Poincar\'e duality (see \cite{Sr} \S 5.6) that $G_i(X/U)
\cong K_i(X/U)$ and $G_i(\overline{X/U}) \cong K_i(\overline{X/U})$.
Since $X/U$ is smooth and $\aone$-contractible, it follows that
$K_i(\Spec k) \longrightarrow K_i(X/U)$ is an isomorphism.

\begin{lem}
\label{lem:short} If $X/U$ is $\aone$-contractible, the localization
sequence gives a short exact sequence
$$
0 \longrightarrow G_0(Z) \longrightarrow G_0(\overline{X/U})
\longrightarrow \Z \longrightarrow 0.
$$
\end{lem}

\begin{proof}
Since $X/U$ is $\aone$-contractible, it follows that $G_1(X/U) \cong
G_1(\Spec k) \cong k^*$, $G_0(X/U) \cong G_0(\Spec k) \cong \Z$.  We
just need to show that the boundary map $G_1(X/U) \longrightarrow
G_0(Z)$ is trivial, or equivalently, that the morphism
$G_1(\overline{X/U}) \longrightarrow G_1(X/U)$ is surjective.  To
see this, observe that each pair $({\mathcal V},\alpha)$ consisting
of a vector bundle on $\overline{X/U}$ and an automorphism $\alpha$
of ${\mathcal V}$ represents an element of $G_1(\overline{X/U})$.
Now, the map $G_1(\overline{X/U}) \longrightarrow G_1(X/U)$ is
induced by restriction.  Since $G_1(X/U) \cong k^*$, we can
represent any class in this group by a pair consisting of a trivial
bundle and an automorphism corresponding to multiplication by an
element of $k^*$.  Such a pair can be extended to give a class in
$G_1(\overline{X/U})$.  (This produces a splitting of the map
$G_1(\overline{X/U}) \longrightarrow G_1(X/U)$ by the canonical
morphism $G_1(\Spec k) \longrightarrow G_1(\overline{X/U})$).
\end{proof}

\begin{lem}
\label{lem:exist} If $X/U$ is a smooth $\aone$-contractible open
dense subscheme of a finite type smooth scheme $\overline{X/U}$,
then there exists a non-trivial vector bundle on $\overline{X/U}$.
\end{lem}

\begin{proof}
First, observe that $G_0(Z)$ is always non-trivial.  Thus, using
Lemma \ref{lem:short}, the map $j^*: K_0(\overline{X/U})
\longrightarrow K_0(X/U)$ always has a kernel.  In particular,
$K_0(\overline{X/U})$ has a generator which is not isomorphic to
$[\O_{\overline{X/U}}]$.  Choosing any vector bundle representing
the isomorphism class of this non-trivial generator gives the
result.
\end{proof}

Now under the hypothesis that the boundary component is of
codimension at least two, non-isomorphic bundles on $\overline{X/U}$
will restrict to non-isomorphic bundles on $X/U$.  Indeed, this
follows by the following ``well-known" result about restrictions of
vector bundles on normal varieties.

\begin{lem}
\label{lem:dense} Assume that the complement of $X/U$ in
$\overline{X/U}$ is of codimension at least two. Then the restriction functor $j^*: {\sf
Vec}(\overline{X/U}) \longrightarrow {\sf Vec}(X/U)$ is
fully-faithful.  Furthermore, the Picard groups of $X/U$ and
$\overline{X/U}$ are isomorphic.
\end{lem}

\begin{proof}
Since $X/U$ is dense in $\overline{X/U}$, the restriction functor is
faithful (any morphism is uniquely determined by restriction to the
generic point).  To check that the functor is full, it suffices to
show that given any pair of locally free sheaves ${\mathcal V}_1$
and ${\mathcal V}_2$ on $\overline{X/U}$, any morphism
$\varphi|_{X/U}: {\mathcal V}_1|_{X/U} \longrightarrow {\mathcal
V}_2|_{X/U}$ extends to a morphism $\varphi: {\mathcal V_1}
\longrightarrow {\mathcal V}_2$.

By assumption, $X/U$ has complement of codimension $\geq 2$ in
$\overline{X/U}$, and $\overline{X/U}$ is smooth and hence normal.
Observe that the canonical morphism $\O_{\overline{X/U}}
\longrightarrow j_*\O_{X/U}$ is an isomorphism (since regular
functions on $X/U$ extend to regular functions on $\overline{X/U}$
by normality). Given ${\mathcal V}_i$ as above, we can choose an
open cover $U_i$ of $\overline{X/U}$ on which ${\mathcal V}_i$
trivialize.  Consider the induced open cover $X/U \cap U_i$ of
$X/U$.  Any morphism $\varphi_{X/U}: {\mathcal V}_1|_{X/U}
\longrightarrow {\mathcal V}_2|_{X/U}$ is specified by a matrix of
regular functions on the $X/U \cap U_i$.  By the extension property
of regular functions mentioned above, this matrix extends uniquely
to give a morphism $\varphi: {\mathcal V}_1|_{U_i} \longrightarrow
{\mathcal V}_2|{U_i}$; furthermore, these morphisms glue to give the
required extension.

Interpreting line bundles in terms of \u Cech cocycles, the
extension property of regular functions on smooth schemes shows that
line bundles on $X/U$ extend to line bundles on $\overline{X/U}$.
\end{proof}

\begin{rem}
Note that any quasi-affine variety admits a canonical open immersion
into the spectrum of its ring of regular functions, which (by
definition of a geometric quotient) here is isomorphic to
$\Spec(k[X]^U)$. Although this last scheme is affine by definition,
it is well-known (Hilbert's 14th Problem) that it need not be
Noetherian, though it is known to be locally of finite type (whence
our conventions). We established in \cite{AD1} Corollary 3.18 (iii)
that in fact the complement of $X/U$ in $\Spec(k[X]^U)$ consists of
codimension $\geq 2$ affine subschemes.  In particular, whenever
$k[X]^U$ is finitely generated, then $X/U$ admits codimension $\geq
2$ affine ``closures".

The smoothness hypothesis on the partial compactification we impose
is for the technical convenience of identifying K-theory with
$G$-theory (via Lemma \ref{lem:exist}), and almost certainly could
be removed.  For the general case, (when $\Spec k[X]^U$ is neither
smooth nor finitely generated), one would have to replace $G$-theory
by Thomason's K-theory (see \cite{Th}). We believe all the lemmas,
with the exception of Lemma \ref{lem:exist}, go through in this
setting: one needs a more subtle argument to extract a vector bundle
from a class in Thomason K-theory.  In any case, the above results
hold for any $\aone$-contractible smooth variety $Y$ that admits a
smooth partial compactification $\overline{Y}$ (we assume
$\overline{Y}$ is a variety) such that $\overline{Y} \setminus Y$
has codimension $\geq 2$ in $\overline{Y}$.
\end{rem}

\subsubsection*{Proof of Theorem \ref{thm:main}}
\begin{proof}[Proof of Theorem \ref{thm:main} (i)]
Suppose $X$ is a smooth affine $\aone$-contractible scheme admitting
a free everywhere stable action of a unipotent group $U$.  Since
unipotent groups in positive characteristic can have non-trivial
finite subgroups, everywhere stability implies properness of the
action of $U$ on $X$ but not necessarily that the action is free.
Imposing this additional condition, the same proof as that of
Theorem 3.10 of \cite{AD1} shows that in this situation a quotient
$X/U$ exists as a quasi-affine smooth scheme (indeed, the proof of
{\em loc. cit.} shows existence of such a quotient is equivalent to
the action being free and everywhere stable).

Now, using \S 3 Example 2.3 of \cite{MV}, together with fact
that $U$ is a special group (i.e., all $U$-torsors are Zariski
locally trivial) we conclude that furthermore $X/U$ is an
$\aone$-contractible smooth scheme (see also \cite{AD1} Key Lemma 3.3).
\end{proof}

\begin{proof}[Proof of Theorem \ref{thm:main} (ii), (iii)]
Statements (ii) and (iii) follow immediately from Lemma
\ref{lem:affinequot}.
\end{proof}

\begin{proof}[Proof of Theorem \ref{thm:main} (iv)]
Combining Lemmas \ref{lem:short}, \ref{lem:exist}, and
\ref{lem:dense}, we obtain the required non-trivial vector bundle on
$X/U$.  We can refine this statement however.  Note that
$Pic(\overline{X/U})$ is necessarily trivial by Lemma
\ref{lem:dense}, thus the non-trivial generator corresponds to a
vector bundle of rank $m \geq 2$.  Furthermore, the vector bundle
representing the non-trivial class on $\overline{X/U}$ is not stably
trivial either so taking direct sums with the trivial bundle
produces non-trivial vector bundles on $\overline{X/U}$ of any rank
$\geq m$.  Restricting these bundles to $X/U$ produces non-trivial
vector bundles of that same rank.
\end{proof}

\begin{rem}
Take $X = {\mathbb A}^n$.  If a unipotent group $U$ acts freely and
everywhere stably on $X$ with a strictly quasi-affine quotient
${\mathbb A}^n/U$ satisfying the hypotheses Theorem \ref{thm:main},
we know there is a non-trivial vector bundle on ${\mathbb A}^n/U$.
The pull-back of this vector bundle to ${\mathbb A}^n$ is
necessarily trivial, thus we see that pull-back by the quotient
morphism does not induce an injection on isomorphism classes of
vector bundles, at least once we have shown such an example exists.
\end{rem}

\section{``Lower bounds" on the failure of $\aone$-invariance}
\label{s:examples}

We must now show that Theorem \ref{thm:main} describes a large class
of schemes.  We are not trying to present a classification, so we
simply supply a class of examples that work over an arbitrary field.

The basic idea is to rewrite the problem of finding $U$-torsors with
$\mathbb{A}^n$ as total space by ``linearizing", i.e., by
restricting from a linear $U$-representation $W$ to a $U$-invariant
subvariety isomorphic to $\mathbb{A}^n$. The simplest cases to
consider are $\ga$-equivariant closed immersions of $\mathbb{A}^n$
as hypersurfaces in $W$, chosen so that $\mathbb{A}^n$ inherits the
structure of a $\ga$-torsor from a larger $\ga$-torsor -- namely an
appropriate open subscheme $W' \subset W$. More specifically, the
geometric points of $W'$ will be ``stable" points of $W$ with
trivial isotropy; such a set can be explicitly identified with the
help of a modified Hilbert-Mumford numerical criterion from
geometric invariant theory (GIT). Furthermore, given such a
$\ga$-equivariant closed immersion, we can identify the ``boundary"
locus (i.e., complement) of $\mathbb{A}^n/\ga$ in
$\Spec(k[\mathbb{A}^n]^{\ga})$ again using geometric invariant
theory.  (For more details on this point of view, we refer the
reader to \cite{AD1} and \cite{DK}.)  With some care one can thereby
arrange that the conditions of Theorem \ref{thm:main} are satisfied.

Let $V$ denote the standard $2$-dimensional representation of
$SL_2$.  By abuse of notation, we will write $V$ instead of
${\mathbb A}(V)$, and furthermore, we choose coordinates $u,v$ on
$V$ throughout and write $0$ for the origin of $V$.  We embed $\ga
\hookrightarrow SL_2$ as the subgroup of {\em lower} triangular
matrices.  Recall that $SL_2/\ga \cong V \setminus 0$; thus $V$ is
an $SL_2$-equivariant completion of $SL_2/\ga$, and we can identify
the identity coset $[e]$ in $SL_2/\ga$ with $\{(0,1)\} \in V$.

Via this embedding $\ga \hookrightarrow SL_2$, an arbitrary
$SL_2$-representation $W$ can be considered as a
$\ga$-representation by restriction. Any given $SL_2$-orbit in $V
\times W$ is contained either in $0 \times W$ or the complement; if
it is in the complement, then it restricts to a $\ga$-orbit in $[e]
\times W$. Similarly any $SL_2$-invariant subscheme of $(V \setminus
0) \times W$ restricts to a $\ga$-invariant subscheme of $[e] \times
W$. We argued in \cite{AD1} (Theorem 3.10, Lemma 4.5, and Theorem
4.11), assuming $k$ was of characteristic $0$, using faithfully flat
descent and the functoriality for quasi-affine maps in GIT, that if
the geometric points of an $SL_2$-invariant subscheme $Y$ in $(V
\setminus 0) \times W$ are stable for the $SL_2$-action on $V \times
W$, then the corresponding $\ga$-invariant subscheme $X$ in $W$ is a
$\ga$-torsor over a quasi-affine variety.  Furthermore we showed the
complement of $X/\ga$ in $\Spec(k[X]^{\ga})$ is the GIT
$SL_2$-quotient of the boundary of the closure $\overline{Y}$ of $Y$
in $V \times W$ (i.e., the complement of the quotient is the
quotient of the complement). Note that $\ga$ acts freely on $X$ if
and only if $SL_2$ acts freely on $Y$. As we explained in the proof
of Theorem \ref{thm:main} (i), if furthermore $Y$ was contained in
the open subscheme of $V \times W$ where $SL_2$ acts freely, the
same result holds in arbitrary characteristic.

\begin{cor}\label{cor:discrete}
For any integers $n \geq 4$, and any $m \geq 1$, there exists a
strictly quasi-affine $\aone$-contractible smooth scheme of
dimension $n$ with at least $m$ non-isomorphic, stably trivial,
non-trival vector bundles of every rank $l$ for sufficiently large
$l$.
\end{cor}

\begin{proof}

Let $W = V^{\oplus 3}$ with coordinates $\{w_1, \ldots, w_6\}$. The
Hilbert-Mumford numerical criterion applied to $SL_2$-orbits in the
$SL_2$-representation $V \times W$ and then restricted to
$\ga$-orbits in $[e] \times W$, implies, as justified above, that
all geometric points in the complement of the subscheme defined by
$\{w_1=0, w_3=0, w_5=0\}$ are stable for the $\ga$-action; in
particular, $V \times W \setminus \{w_1=0, w_3=0, w_5=0\}/\ga$ is a
quasi-affine geometric quotient.

For $f + 1$ a $1$-variable polynomial with no repeated roots and
constant term $1$, consider a $\ga$-invariant hypersurface $X$ given
by the $\ga$-invariant equation $w_1 = 1 + f(w_3 w_6-w_4 w_5)$.
Observe first that by the preceding paragraph all geometric points
of $X$ are stable, since any non-stable point must satisfy
$w_1=w_3=w_5=0$ which violates the hypersurface equation because $0
\neq 1$.  It follows that the geometric quasi-affine quotient
$X/\ga$ exists. Second, $X$ is isomorphic to $\mathbb{A}^5$ via the
closed immersion $w_2 = z_1, \ldots, w_6 = z_5$, where $\{z_1,
\ldots, z_5\}$ are coordinates on $\mathbb{A}^5$. Third, note that
$f(w_3 w_6-w_4 w_5)$ is $SL_2$-invariant; it follows easily that the
associated hypersurface equation defining $\overline{Y}$ in $V
\times W$ is $uw_2-vw_1 = 1 + f(w_3 w_6-w_4 w_5)$, where $(u,v)$ are
the coordinates on the first factor of $V$. Fourth, again by the
Hilbert-Mumford criterion, all geometric points of $\overline{Y}$
are stable with respect to the $SL_2$-action on $V \times W$, so in
particular the $SL_2$ action on $\overline{Y}$ is proper, and the
affine geometric quotient $\overline{Y}/SL_2$ exists.

The boundary $B = \overline{Y} \setminus Y$ of $\overline{Y}$,
equivalently the closed subscheme of $\bar{Y}$ defined by the
simultaneous vanishing of $u$ and $v$, is explicitly given by $f(w_3
w_6-w_4 w_5)+1=0$ in $0 \times W$; it is $SL_2$-invariant,
codimension $2$ in $\overline{Y}$, and its geometric points are all
$SL_2$-stable. As per the discussion preceding this Corollary, this
means the complement of $X/\ga = \mathbb{A}^5/\ga$ in
$\Spec(k[\mathbb{A}^5]^{\ga})$ is given by $B/SL_2$, which is
necessarily again codimension $2$.

Observe that by the Jacobian criterion, the fact that $f+1$ has no
repeated roots implies that both $\overline{Y}$ and $B$ are smooth
schemes for any $k$.  Also, $B$ has $m$ disjoint components, where
$m$ is the degree of $f$.  Recall that an action is called
set-theoretically free if its stabilizers at $k$-points are trivial.
Furthermore, proper, set-theoretically free actions are free ({\em
cf.} \cite{AD1} Lemma 3.11).  A direct computation shows the only
geometric points having non-trivial stabilizers for the $\ga$-action
on $W$ lie in the non-stable locus, so $\ga$ acts freely on $X$ and
hence $SL_2$ acts freely on $Y$. Also note that $SL_2$ acts freely
on $B$ in $0 \times W$; indeed $B$ is a finite disjoint union of
rank $2$ vector bundles over $SL_2$.  Consequently $SL_2$ acts not
only properly but also set-theoretically freely, hence freely, on
all of $\overline{Y}$, so $\overline{Y}/SL_2$ is smooth.  Denote the
$SL_2$-quotient of the boundary $B$ by $Z$; then $Z$ is smooth and
codimension 2 in $\overline{Y}/SL_2 \cong
\Spec(k[\mathbb{A}^5]^{\ga})$.

Indeed, $Z$ in these examples is isomorphic to a disjoint union of
$m$ copies of the affine plane.  It follows that for any $m \geq 1$,
we can choose $f$ and hence $X$ so that $K_0(Z)$ is isomorphic to
$\Z^{\oplus m}$. Thus Lemma \ref{lem:short} shows that $K_0(\Spec
k[\mathbb{A}^5]^{\ga}) \cong \Z^{\oplus m +1}$.  As vector bundles
representing different classes in $K_0$ are not stably equivalent,
by taking direct sums with trivial bundles we get $m$ non-trivial,
non-isomorphic bundles in every sufficiently large rank.  Then
restriction, by Lemma \ref{lem:dense}, gives the desired bundles on
$\mathbb{A}^5/\ga$.

Higher dimensional examples immediately follow by taking other
representations; for example, $W = V^{\oplus 3} \oplus k^r$, where
$k$ denotes the trivial representation and $\mathbb{A}^{5+r}$ is
presented as a hypersurface with the same equation as above.
\end{proof}

\begin{rem}
We expect it is possible to construct a smooth quasi-affine
$3$-dimensional variety with the desired properties via unipotent
quotients of an affine space. However, we do not believe there are
any smooth $\aone$-contractible surfaces other than $\mathbb{A}^2$
(see also Remark \ref{rem:dim2}); this is known to be true over $\cplx$ (see \cite{AD2}).
\end{rem}

\begin{rem}
The quasi-affine quotient scheme in the simplest case of the
construction from Corollary \ref{cor:discrete} (where $f$ is the
identity so that $\mathbb{A}^5$ is defined by $w_1 = 1 +
(w_3w_6-w_4w_5)$) is very pleasant to visualize.  An easy
computation with invariants presents $Spec(k[X]^{\ga})$ as a quadric
hypersurface in $\mathbb{A}^5$. When $k=\mathbb{C}$ this may be
thought of as the complexification of a sphere, that is, $T^*(S^4)$.
Here $B$ is a single affine plane: over $\mathbb{C}$, $B$ is the
cotangent plane at a point, so the complement of $B$ is clearly
contractible as a complex manifold. This particular example of a
quasi-affine contractible complex variety, with a different
presentation, was known to Winkelmann \cite{Wi}.
\end{rem}

In the other direction, given a desired boundary we can often pick
the defining hypersurface equation for $X$ so as to yield a quotient
with that specified boundary. Varying the boundary in a family may
be realized by varying the defining hypersurface equation in the
fixed $\ga$-representation $W$.  By arranging for a boundary $Z$
with a large $K_0(Z)$, we can then by the above process get smooth
$\mathbb{A}^1$-contractible schemes with arbitrarily many
non-isomorphic vector bundles, and indeed find arbitrary dimensional
families of such schemes.

\begin{cor}\label{cor:moduli}
For any integers $n \geq 6$, $m \geq 1$, and $l \geq 1$ there exists
an $m$-dimensional smooth scheme $S$ and a smooth morphism $f: X
\longrightarrow S$ of relative dimension $n$ whose fibers are
strictly quasi-affine $\aone$-contractible smooth schemes, pair-wise
non-isomorphic, each of which possesses at least $l$-dimensional
moduli of stably trivial, non-trivial vector bundles in every
suitably large rank.
\end{cor}

\begin{proof}
We use notation and terminology as in the proof of Corollary
\ref{cor:discrete}.  Consider $W = V^{\oplus 4}$, with coordinates
$\{w_1, \ldots, w_8\}$, as an $SL_2$-representation and hence a
$\ga$-representation (where $\ga \hookrightarrow SL_2$ as lower
triangular matrices, as before). Then the hypersurface $w_1 = 1 +
f(w_3w_6-w_4w_5, w_3w_8-w_4w_7,w_5w_8-w_6w_7)$ in $W$ is isomorphic
to $\mathbb{A}^7$; the closed immersion is determined by function
$w_2 = z_1, \ldots, w_8=z_7$, where $\{z_1, \ldots, z_7\}$ are the
coordinates on $\mathbb{A}^n$. It is easily checked that the
restriction of the linear $\ga$-action on $W$ to this $\mathbb{A}^7$
hypersurface is everywhere stable, by using $SL_2$-stability for $V
\times W$ as before.

Note that $f(w_3w_6-w_4w_5, w_3w_8-w_4w_7,w_5w_8-w_6w_7)$ is an
$SL_2$-invariant, so the associated hypersurface equation defining
$\overline{Y}$ in $V \times W$ is $uw_2-vw_1 = 1+ f(w_3w_6-w_4w_5,
w_3w_8-w_4w_7,w_5w_8-w_6w_7)$. Over any field $k$, for generic $f$
this describes a smooth hypersurface in the $SL_2$-stable locus of
$V \times W$, all of whose points have trivial isotropy in $SL_2$;
we leave the details to the reader.  In particular for generic $f$
the $SL_2$ action on $\overline{Y}$ is free, and the quotient
$\overline{Y}/SL_2$ is smooth.  Since the boundary $B$ is defined by
the simultaneous vanishing of $u$ and $v$, it is a hypersurface in
$0 \times W$ and so is codimension $2$ in $\overline{Y}$. Indeed,
$B$ consists of a rank $2$ vector bundle over a principal
$SL_2$-bundle over a smooth affine surface. The quotient $Z =
B/SL_2$ is thus codimension $2$ and a smooth subvariety of the
smooth $\overline{Y}/SL_2$.  Consequently if $Y_1$ and $Y_2$
(respectively, $Z_1$ and $Z_2$) are the $SL_2$-invariant varieties
(respectively, boundaries of the quotients) associated with two
different choices of $f$, say $f_1$ and $f_2$, then any morphism
from $Y_1/SL_2$ to $Y_2/SL_2$ extends to a morphism from
$\overline{Y_1}/SL_2$ to $\overline{Y_2}/SL_2$ and vice-versa; so
$Y_1/SL_2 \cong Y_2/SL_2 \Rightarrow \overline{Y_1}/SL_2 \cong
\overline{Y_2}/SL_2 \Rightarrow Z_1 \cong Z_2$.  In particular, if
$Z_1 \not\cong Z_2$ then $Y_1/SL_2 \not\cong Y_2/SL_2$. Thus the
fact that there are arbitrary dimensional moduli of the surfaces
$1+f(x,y,z)=0$, and hence of the boundaries $Z$, means there are
arbitrary dimensional moduli of $Y/SL_2 \cong \mathbb{A}^7/\ga$
associated with varying the $\ga$-action ({\em cf.} \cite{AD1} Lemma
5.5)).

Since $Z \cong B/SL_2$ is a vector bundle over a smooth affine
surface $S$, the map $K_0(Z) \longrightarrow K_0(S)$ is an
isomorphism. Furthermore, the smooth affine surface is defined as a
hypersurface in ${\mathbb A}^3$.  So for example, if we take a
hypersurface isomorphic to a product of a smooth affine curve and
the affine line (the reader may check that a family of examples in
any genus may be chosen so that $\overline{Y}$ is smooth and so that
$B$ is contained in the open subscheme of $0 \times W$ on which
$SL_2$ acts freely, thus guaranteeing $Z$ is smooth), we see that
$K_0(Z)$ can be made arbitrarily large by making the genus of the
curve high: specifically, line bundles are cancellation stable so
there is an injection from $Pic(Z)$ into $K_0(Z)$, and affine curves
have moduli of line bundles of dimension increasing with the genus.
Because everything is smooth, the same argument as in the previous
Corollary now implies the desired statement for $\ga$-quotients of
$\mathbb{A}^7$; quotients for larger dimensional $\mathbb{A}^n$ may
be achieved by taking other representations, e.g., $W \oplus k^r$
for $k$ the trivial representation and the defining equation the
same as above.
\end{proof}

\section{Some comments on the generalized Serre problem}
\begin{prop}
Suppose $X$ is a topologically contractible smooth complex variety
of dimension $\leq 2$, then every vector bundle on $X$ is isomorphic
to a trivial bundle.\footnote{Added in proof: Proposotion \ref{prop:GS} can be found in Corollary 2 of \cite{GS1} with a similar proof.}
\label{prop:GS}
\end{prop}

\begin{proof}
If $X$ is a topologically contractible smooth complex curve, then
$X$ is isomorphic to the affine line and the result follows from the
Quillen-Suslin theorem.  Therefore, we can assume that $X$ has
dimension $2$.  Suppose therefore that $X$ is a topologically
contractible smooth complex surface.

By a Lemma of Fujita (see e.g., \cite{Za} Lemma 2.1), we know that
any such surface is affine.  By a Theorem of Gurjar-Shastri, (see
\cite{Za} Theorem 2.1) we know that any topologically contractible
smooth complex surface is rational.  In particular, $X$ admits a
smooth projective compactification $\bar{X}$ which is a smooth
projective rational surface.  By the classification of surfaces
$\bar{X}$ is birationally equivalent to a ruled surface. Murthy (see
\cite{Murthy} Theorem 3.2) has shown that every vector bundle on any
affine surface birationally equivalent to a ruled surface is
necessarily isomorphic to the direct sum of a trivial bundle and a
line bundle.  Thus, if $Pic(X)$ is trivial, it follows that every
vector bundle on $X$ is isomorphic to a trivial bundle.

To see that $Pic(X)$ is trivial for a smooth contractible surface,
choose a compactification of $\bar{X}$ whose boundary is a simple
normal crossings divisor $D$.  We have an exact sequence for Chow
groups
$$
CH_i(D) \longrightarrow CH_i(\bar{X}) \longrightarrow CH_i(X) \longrightarrow 0.
$$
In particular, taking $i = 1$ and using the fact that $\bar{X}$ and
$X$ are smooth, we see that $Pic(\bar{X}) \longrightarrow Pic(X)$ is
surjective.  By Corollary 2.2 of \cite{Za}, we know that
$Pic(\bar{X})$ is freely generated by the irreducible components of
$D$ and thus $Pic(X)$ is necessarily trivial.
\end{proof}

\begin{cor}
\label{cor:strictqafftrivvb} If $X$ is any topologically
contractible smooth complex algebraic surface, and $p_1,\ldots,p_n$
are finitely many points on $X$, then all vector bundles on $X
\setminus \setof{p_1,\ldots,p_n}$ are trivial.
\end{cor}

\begin{proof}
Indeed, if $\F$ is a locally free sheaf on $X \setminus
\setof{p_1,\ldots,p_n}$, then there always exists a coherent
extension $\bar{\F}$ of $\F$ to $X$.  The double dual
$\bar{\F}^{\vee\vee}$ is a reflexive sheaf on $X$, which must be
locally free since $X$ is a smooth surface; this provides a locally
free extension of $\F$.  We have just shown that all vector bundles
on such an $X$ are in fact trivial, and thus $\F$ must be a trivial
bundle as well.
\end{proof}

\begin{cor}\label{cor:contrmoduli}
There are positive dimensional moduli of smooth algebraic surfaces
which admit only trivial vector bundles.  These can be chosen so
that they are affine and, as complex manifolds, contractible or
quasi-affine and non-contractible.
\end{cor}

\begin{proof}
There are contractible smooth affine algebraic surfaces of log
Kodaira dimension 1 that admit deformations (see \cite{FlZa} Example
6.9).  Upon removing finitely many points, Corollary
\ref{cor:strictqafftrivvb} finishes the result.
\end{proof}

\begin{rem}
\label{rem:dim2} None of the examples mentioned in the proof of the
Corollary are $\aone$-contractible.  In fact, contractible smooth
surfaces of positive log Kodaira dimension are known not to be
$\aone$-contractible (see \cite{AD2}).  Roughly speaking, this is
the case because positive log Kodaira dimension surfaces do not have
``many" rational curves; in order for a variety to even be
$\aone$-connected, one expects that it should be covered by chains
of $\aone$s.
\end{rem}

\begin{rem}
For topologically contractible smooth complex affine varieties of
dimension $n \geq 3$, the functor $X \mapsto {\mathscr V}(X)$
becomes even more subtle.  Results of Suslin imply that for
projective modules of rank $\geq n$ stable isomorphism implies
isomorphism.  In this direction, results of Murthy imply(see
\cite{Murthy2} Corollary 2.11) that if $f,g$ are elements of the
polynomial ring $\cplx[x_1,\ldots,x_n]$ with $g \neq 0$, then all
stably free modules over $\cplx[x_1,\ldots,x_n,f/g]$ of rank $\geq
n-1$ are free.  If $X$ is a topologically contractible smooth
complex $3$-fold, of this form, then all vector bundles on $X$ are
trivial if and only if $Pic(X)$ is trivial.  In particular, Murthy
(\cite{Murthy2} Theorem 3.6) uses this to deduce that all the
Koras-Russell threefolds, in particular the famous Russell cubic
surface $x + x^2y + z^2 + t^3 = 0$, satisfy the generalized Serre
problem.  At the moment, it is not known whether or not the Russell cubic
is $\aone$-contractible.
\end{rem}

\begin{footnotesize}
\bibliographystyle{alpha}
\bibliography{exoticaffinesIII}
\end{footnotesize}

\end{document}